\documentclass[12pt]{article}
\usepackage{amsmath,amsfonts}
\begin{document}

\baselineskip 20pt
\date{ December 29, 2014}
\title{On the zero-free polynomial approximation problem}

\author{Arthur A.~Danielyan}

\maketitle

\begin{abstract}
\noindent Let $E$ be a compact set in $\mathbb C$ with connected complement, and let $A(E)$ be the class
of all complex continuous function on $E$ that are analytic in the interior $E^0$ of $E$. Let $f \in A(E)$ be zero free on $E^0$. 
By  Mergelyan's theorem $f$ can be uniformly approximated on $E$ 
by polynomials, but is it possible to realize such approximation by polynomials that are zero-free on $E$?
This natural question has been proposed by J. Andersson and  P. Gauthier. 
So far it has been settled for some particular sets $E$.
The present paper describes classes of functions for which zero free approximation is possible on an arbitrary $E$.
\end{abstract}

\begin{section}{Introduction and the main result.}

For a set $M$ in the complex plane $\mathbb C$ we denote as usual by  $M^0$, $\partial M$, and $\overline M$
the interior, the boundary, and the closure  of $M$, respectively.

Suppose $E$ is an arbitrary compact subset in $\mathbb C$ such that $\mathbb C \setminus E$ is connected.
Let $A(E)$ be the usual space (algebra) of all complex-valued continuous functions on $E$ that are analytic in $E^0$.
The following well-known approximation theorem is due to S.N. Mergelyan (see \cite{Mer} or \cite{Rud}).

 \vspace{0.25 cm}

{\bf Theorem A.} {\it Let 
$f \in A(E)$. Then for each $\epsilon >0$ there exists a polynomial $P(z)$ such that $|f(z) - P(z)|< \epsilon$ for all $z \in E$. }

 \vspace{0.25 cm} 

 In regards to this theorem, the following
 natural question is on the possibility of approximation by polynomials that are zero-free on the set of approximation.
 
 \vspace{0.25 cm} 
 
 {\bf Question 1.} Let $f \in A(E)$ has no zeros on $E^0$ and let $\epsilon >0$.
 Does there exist a polynomial $P(z)$ with no zeros on $E$ such that $|f(z) - P(z)|< \epsilon$ for all $z \in E$? 
 
 \vspace{0.25 cm} 
 
 
 Note that if $f \in A(E)$  has a zero at a point of $E^0$, then by Hurwitz's theorem for such $f$ a zero free polynomial approximation is impossible.
 
Question 1 has been proposed by 
 J. Andersson and P. M. Gauthier (cf. \cite{And2}) and it has been investigated in the recent papers \cite{And1}, \cite{And2}, \cite{AndGa}, \cite{Gaut1},
 \cite{Gaut2},  \cite{GaKn}, \cite{Khr}.
 Several interesting results of affirmative character have been proved under various restrictions on $E$, but the question still remains open in the general case. 
 

The present paper describes classes of functions defined on an arbitrary compact set $E$ with connected complement 
for which zero free approximation is possible. 

To avoid any possible confusion, let us first remind (introduce) some terminology.

\vspace{0.25 cm}

{\bf Definition 1.} Let $E \subset \mathbb C$ be a compact set with connected complement.
We call $H \subset \partial E$ a zero set if  some $g \in A(E)$
vanishes precisely on $H$ (that is, $H=\{z \in E: g(z)=0\}$).

\vspace{0.25 cm}

Clearly any such $H$ is a closed subset on $\partial E$ and also the function $g$ is zero free on $E^0$.

The main result of this paper is the following:

\vspace{0.25 cm}

{\bf Theorem 1}. {\it Let $E \subset \mathbb C$ be a compact set with connected complement and let $H \subset \partial E$ be a zero set. 
Then there exists $f \in A(E)$ with $H$ as its zero set and allowing uniform approximation by zero free on $E$ polynomials.}

  \vspace{0.25 cm}
 
 The theorem implies that for any prescribed zero set $H$ the zero free approximation is possible at least for some functions.
In particular, a possible counterexample to Question 1, except the set $H$, has to depend also on a specific function of $A(E)$ vanishing on $H$.
(Since in the general case both $E$ and $H$ may be complicated sets, the last conclusion sheds some light on the nature of a possible 
counterexample; see also the discussion in \cite{And2}, Sect. 4.)

Each polynomial $P(z) = a(z-\zeta_1)(z- \zeta_2)...(z-\zeta_n)$ can be considered as a continuous function of its zeros $\zeta_1, \zeta_2, ..., \zeta_n$ (and $z$).
This, combined with the fact that $\partial E$ is nowhere dense in $\mathbb C$, clearly implies that Question 1 is equivalent to the following question.

\vspace{0.25 cm} 
 
 {\bf Question 1'.} Let $f \in A(E)$ has no zeros on $E^0$ and let $\epsilon >0$.
 Does there exist a polynomial $P(z)$ with no zeros on $E^0$ such that $|f(z) - P(z)|< \epsilon$ if $z \in E$? 
 
 \vspace{0.25 cm}


Assume $E^0=\emptyset$ and $f \in A(E)$. Then by Theorem A (or merely by its particular case known as 
Lavrentiev's theorem) a polynomial approximation to $f$ is possible, and trivially Question 1' has an affirmative answer (since $E^0=\emptyset$).
Then the same is true for Question 1 and thus we have the following result of \cite{And1} (cf. also \cite{And2}): 

 \vspace{0.25 cm} 
 
{\bf Proposition 1.} {\it If $E^0=\emptyset$, zero free polynomial approximation on $E$ is always possible.}

 \vspace{0.25 cm} 
 
If  $E^0 \neq \emptyset$ but $f \in A(E)$ is zero free on $\overline {E^0}$ then obviously by Theorem A one can approximate $f$ on $E$ by a 
polynomial which is zero free on $\overline {E^0}$.
This implies, in particular, that Question 1' has  an affirmative answer. Then, as above, Question 1 too has an affirmative answer and we arrive to the following extension of Proposition 1:
 
  \vspace{0.25 cm} 

{\bf Proposition 2}.  {\it If $f \in A(E)$ has no zeros on $\overline {E^0}$ then $f$ can be uniformly approximated on $E$ by zero 
free (on $E$) polynomials.}\footnote{Note that Proposition 2 immediately implies the second part of the proof of Theorem 5 of \cite{And2} 
since in \cite{And2}, p. 206, the (extended) function $H$ on $K$ satisfies the conditions of Proposition 2.}
 
 \vspace{0.25 cm}

The following simple ingredient of the proof of Theorem 1, in fact, is also
 a sufficient condition for zero free approximation on an arbitrary compact set with connected complement.

 \vspace{0.25 cm}
 
{\bf Theorem 2}.  {\it Let $E$ be as in Theorem 1. If $f =u+iv \in A(E)$ is zero free on $E^0$ and $u$ (or $v$) is  
either nonnegative on $E^0$ or nonpositive on $E^0$,
then $f$ allows uniform approximation by zero free on $E$ polynomials.}

 \vspace{0.25 cm}

 We use the following well-known peak-interpolation theorem of E. Bishop (cf. \cite{Rud}, p. 135).

 \vspace{0.25 cm} 

{\bf Theorem B.} {\it Let $X$ be a compact Hausdorff space and let $A$ be a closed linear subspace 
of the space $C(X)$ of all complex continuous functions on $X$. Suppose $K$ is a compact $G_\delta$ subset of $X$ such that
$|\mu|(K)=0$ for every regular complex Borel measure on $X$ which is orthogonal to $A$. Then $K$ is a peak-interpolation set for $A$. }

 \vspace{0.25 cm} 
 
 Recall that, by the definition, $K$ is a peak-interpolation set for $A$ if for any continuous $f$ on $K$, not identically zero, there exists $u\in A$ such that $u=f$ on $K$ and $|u(z)|<||f||_K$ for every  $z \in X\setminus K$. Note that a new approach to Theorem B has been presented in \cite{Dan}.
 
 Except the theorems of Mergelyan and Bishop, below we use the classical decomposition formula of orthogonal complex measures as it has been presented  by L. Carleson in his potential theory based well-known research \cite{Carl} on Mergelyan's theorem.

\end{section}

\begin{section}{Proofs}
{\it Proof of Theorem 1.} 
Let $H \subset \partial E$ be a zero set, as stated in the theorem, and let $g \in A(E)$ be a function vanishing precisely  on $H$. Denote by  $D_1, D_2, ... ,D_n, ... $ the connected components of $E^0$. Because the complement of $E$ is connected, each $D_n$ is simply connected. Let $H_n=H\cap \partial D_n$ be the portion of the zero set $H$ on the boundary of $D_n$.  Denote by $\eta_n$ any harmonic measure  of domain $D_n$. Then  $\eta_n(H_n)=0$, because otherwise the function $g$ would be identically zero in $D_n$ (contrary to the hypothesis of the theorem). 

Denote by $B$ the set of restrictions to $\partial E$ of all elements of $A(E)$. By the modulus maximum principle, $B$ is a closed subspace of all continuous functions on $\partial E$.

Consider an arbitrary regular complex Borel measure $\mu $ on $\partial E$, which is orthogonal to $B$. Then $\mu$ in particular is orthogonal to all polynomials. 
By a classical description of such $\mu$,
for each natural $n$ there exists an orthogonal to $B$
measure $\mu_n$ concentrated on $\partial D_n$, which is absolutely continuous with respect to $\eta_n$ and such that $\mu = \sum \mu_n$, the series being convergent in the sense of total variations (cf. \cite{Carl}, p. 175).

By an elementary property, the total variation $|\mu_n|$ too is absolutely continuous with respect to $\eta_n$.
Since $\eta_n(H_n)=0$, we have $|\mu_n| (H_n)=0$. Like $\mu_n$, the total variation $|\mu_n|$ is concentrated on $\partial D_n$,
and therefore $|\mu_n|(H)=0$. This implies that the measure  $\mu_n$ takes the value zero on the subsets of $H$ and the decomposition
$\mu = \sum \mu_n$ shows that
  $\mu$ possesses the same property. Thus, by the definition of $|\mu|$, we have $|\mu|(H)=0$.
 
In Bishop's theorem take $X=\partial E$, $A=B$, and $K=H$. As we have just seen, for $H$ the relation  $|\mu|(H)=0$ for all appropriate measures $\mu$ is satisfied. Since in addition $H$ is a compact (and also a $G_\delta$) subset of $\partial E$, it satisfies the conditions of 
Bishop's theorem. By that theorem, in particular there exists a function $g_1 \in B$ which peaks on $H$, that is, $g_1$ equals to $1$ on $H$ and $|g_1|$ is less than $1$ on $\partial E \setminus H$. Of course, $g_1$ can be extended also on the set $E^0$ to make a function belonging to $A(E)$,
and we have  for this extended function (denoted again $g_1$) that $|g_1|$ is less than $1$ on $E \setminus H$. Let $f = 1 - g_1$. Obviously $f \in A(E)$ and $H$ is the zero set of $f$.

The real part of $f$ is nonnegative on $E$ and for any $\delta >0$ the real part of $\delta +f$ is positive on $E$. The function $\delta +f$ is zero free on $E$ and therefore  $|\delta + f|$ is bounded away of zero on $E$. Now Mergelyan's theorem immediately implies that $\delta + f$ can be uniformly approximated (with any given preciseness) by a polynomial which is  zero free on $E$. Since $\delta>0$ can be taken small enough, the same polynomial will approximate the function $f$ as well. The proof is over.

\vspace{0.25 cm}

{\it Proof of Theorem 2.} Without loss of generality we assume that the real part $u$ of $f$ is nonnegative on $E^0$; then $u$ is nonnegative also on $\overline {E^0}$.
The argument is almost the same as that
 in the last paragraph of the above proof of Theorem 1. Even now the real part of the function $\delta +f$ is positive only on $\overline {E^0}$, things remain simple.
 The shortest way to complete the proof is to apply Proposition 2 for the function $\delta +f$ which is zero free on $\overline {E^0}$.  Again, since $\delta>0$  can be taken arbitrarily small, the (zero free on $E$) polynomial approximating $\delta +f$ also approximates $f$. The proof of Theorem 2 is over.


\end{section}

\begin{minipage}[t]{6.5cm}
Arthur A. Danielyan\\
Department of Mathematics and Statistics\\
University of South Florida\\
Tampa, Florida 33620\\
USA\\
{\small e-mail: adaniely@usf.edu}
\end{minipage}

\end{document}